\documentclass{article}
\usepackage[utf8]{inputenc}
\usepackage{amsmath}
\usepackage{amsfonts} 
\usepackage{amsthm}
\usepackage{verbatim}
\usepackage{etoolbox}
\usepackage{dirtytalk}
\usepackage{tikz-cd}
\usepackage{derivative}
\usepackage{amssymb}
\usepackage[maxnames=99]{biblatex} 
\theoremstyle{plain}
\newtheorem{thm}{Theorem}[section]
\newtheorem{lem}[thm]{Lemma}
\newtheorem{obs}[thm]{Observation}
\newtheorem{prop}[thm]{Proposition}
\newtheorem{cor}[thm]{Corollary}
\newtheorem{claim}{Claim}

\newtheorem*{claim*}{Claim}
\newtheorem{que}[thm]{Question}

\theoremstyle{definition}
\newtheorem{defi}[thm]{Definition}
\newtheorem{defiandobs}[thm]{Definition and Observation}

\theoremstyle{remark}

\AtBeginEnvironment{thm}{\setcounter{equation}{0}}
\AtBeginEnvironment{lem}{\setcounter{equation}{0}}
\AtBeginEnvironment{obs}{\setcounter{equation}{0}}
\AtBeginEnvironment{defi}{\setcounter{equation}{0}}
\AtBeginEnvironment{prop}{\setcounter{equation}{0}}
\AtBeginEnvironment{cor}{\setcounter{equation}{0}}

\AtBeginEnvironment{thm}{\setcounter{case}{0}}
\AtBeginEnvironment{lem}{\setcounter{case}{0}}
\AtBeginEnvironment{obs}{\setcounter{case}{0}}
\AtBeginEnvironment{defi}{\setcounter{case}{0}}
\AtBeginEnvironment{prop}{\setcounter{case}{0}}
\AtBeginEnvironment{cor}{\setcounter{case}{0}}

\AtBeginEnvironment{thm}{\setcounter{claim}{0}}
\AtBeginEnvironment{lem}{\setcounter{claim}{0}}
\AtBeginEnvironment{obs}{\setcounter{claim}{0}}
\AtBeginEnvironment{defi}{\setcounter{claim}{0}}
\AtBeginEnvironment{prop}{\setcounter{claim}{0}}
\AtBeginEnvironment{cor}{\setcounter{claim}{0}}

\newcommand{\C}{\mathbb{C}}
\newcommand{\OO}{\mathcal{O}}

\newcommand{\Z}{\mathcal{Z}}
\newcommand{\Bl}{\textnormal{Bl}}

\title{Uniform rationality of products with the plane}
\author{Juliusz Banecki}
\date{}
\AtEndDocument{%
  \par
  \textsc{Juliusz Banecki,
Faculty of Mathematics and Computer Science,
Jagiellonian University,
ul. Lojasiewicza 6, 30-348 Krakow, Poland}\\
\textit{E-mail address}: 
\texttt{juliusz.banecki@student.uj.edu.pl}
}

\begin{document}
\maketitle
\begin{abstract}
Let $X$ be a nonsingular rational variety. We prove that $X\times \mathbb{C}^2$ is uniformly rational. It follows that nonsingular stably rational varieties are stably uniformly rational.
\end{abstract}
\section{Introduction}
Throughout the paper, the term \emph{variety} refers to a quasi-projective variety over $\C$, i.e. an irreducible Zariski locally closed subset of a projective space $\mathbb{P}_\C^n$. Morphisms of such varieties are called \emph{regular mappings}. 

Let $X$ be a variety. Given a point $x_0\in X$, by $\OO_{X,x_0}$ we denote the ring of germs of complex-valued regular functions at $x_0$, while by $\OO(X)$ we denote the ring of globally defined regular functions on $X$. By a subvariety of $X$ we mean an irreducible Zariski locally closed subset of it.

We are concerned with \emph{uniformity} problems in birational geometry. To specify what we mean by that, we need to introduce a few standard definitions:
\begin{defi}
Let $X$ be a variety. 

We say that $X$ is \emph{rational} if it is birationally equivalent to the projective space $\C^{\dim X}$. Alternatively, $X$ is rational if there exists a nonempty Zariski open subset $U\subset X$ which is biregularly isomorphic to a Zariski open subset $V\subset \C^{\dim X}$.

We say that $X$ is \emph{stably rational} if $X\times \C^n$ is rational for some $n\geq 0$.

Finally, we say that $X$ is \emph{retract rational}, if there exists a natural number $n$, a nonempty Zariski open subset $U\subset X$, a nonempty Zariski open subset $V\subset \C^n$ and two regular mappings $i:U\rightarrow V$ and $r:V\rightarrow U$ such that their composition $r\circ i $ is equal to the identity $\textnormal{id}_U$.
\end{defi}

Each of the above conditions admits the corresponding uniform version:
\begin{defi}
Let $X$ be a variety.

We say that $X$ is \emph{uniformly rational} if for every point $x\in X$ there exists a Zariski open neighbourhood $U_x\subset X$ of $x$ which is biregularly isomorphic to a Zariski open subset $V_x\subset \C^{\dim X}$.

We say that $X$ is \emph{stably uniformly rational} if $X\times \C^n$ is uniformly rational for some $n\geq 0$.

Finally, we say that $X$ is \emph{uniformly retract rational}, if for every $x\in X$, there exists a natural number $n_x$, a nonempty Zariski open neighbourhood $U_x\subset X$ of $x$, a nonempty Zariski open subset $V_x\subset \C^{n_x}$ and two regular mappings $i_x:U_x\rightarrow V_x$ and $r_x:V_x\rightarrow U_x$ such that their composition $r_x\circ i_x $ is equal to the identity $\textnormal{id}_{U_x}$.
\end{defi}

There are the following obvious implications between the different notions:
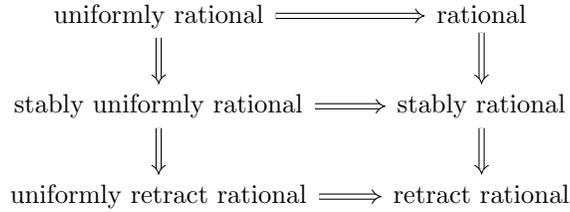
\begin{center}
\begin{tikzcd}
\text{uniformly rational} \arrow[r, Rightarrow] \arrow[d, Rightarrow]        & \text{rational} \arrow[d, Rightarrow]        \\
\text{stably uniformly rational} \arrow[d, Rightarrow] \arrow[r, Rightarrow] & \text{stably rational} \arrow[d, Rightarrow] \\
\text{uniformly retract rational} \arrow[r, Rightarrow]                      & \text{retract rational}                     
\end{tikzcd}
\end{center}

One can naturally pose the following question: which of the implications in the diagram are strict?

When it comes to the vertical ones, it is well known that there exist stably rational varieties which are not rational \cite{BeauvilleColliotTheleneSansucSwinnertonDyer1985}. Passing to an open subset we see that we get examples of stably uniformly rational varieties, which are not uniformly rational.

The case of the other two vertical implications is more subtle; to our best knowledge it is an open problem whether there exists a retract rational variety over $\C$ which is not stably rational, although such examples are known to exist over certain non-closed fields (see \cite[p. 2]{LangeSchreiederLowDegHypersurfaces2024}). 

In this paper, we are concerned with the horizontal arrows. It is easy to see that uniformly retract rational varieties must be nonsingular, so to invert the horizontal arrows we must furthermore assume that $X$ is nonsingular. Other than that, there are no obvious obstructions, so we can ask the following question:
\begin{que}\label{main_que}
Which of the rationality notions (rational, stably rational, retract rational), imply the corresponding uniform notion under the assumption that $X$ is nonsingular?
\end{que}

Of the three questions, so far only the one regarding retract rationality has been answered, by the author in \cite{BaneckiUniformRetract}:
\begin{thm}
Let $X$ be a nonsingular retract rational variety. Then, it is uniformly retract rational.
\end{thm}

The question regarding rational varieties was first posed by Gromov \cite[Sec. 3.5.E\textquotesingle\textquotesingle\textquotesingle]{GromovOka1989} and has been studied since then, although to our best knowledge it is still open (see \cite[Question 1.1]{BogomolovBohning2014} and \cite[Question 4.9]{ZaidenbergAlgebraicGromovEllipticity2024} for a survey). It is worth noting that here that the term \emph{uniformly rational} is not entirely standard; such varieties are called \emph{regular} in \cite{GromovOka1989}, \emph{plain} in \cite{BodnarHauserSchichoVillamayor2008} or \emph{locally flattenable} in \cite{Popov2020VariationsOnTheThemeOfZariskiCancellationProblem}. 

The most promising result in the positive direction regarding the question is the following observation, which first appeared in the paper \cite{GromovOka1989}, where Gromov attributes it to Bogomolov:
\begin{obs}[{\cite[Proposition 3.5.E]{GromovOka1989},\cite[Proposition 2.6]{BogomolovBohning2014}}]\label{obs_Bogomolov}
Let $X$ be a uniformly rational variety, and let $Y\subset X$ be its nonsingular Zariski closed subvariety. Then, the blowup $\Bl_YX$ of $X$ along $Y$ is uniformly rational.
\end{obs}
Conversely, if one could prove the inverse implication (that uniform rationality of $\Bl_Y X$ implies that $X$ is uniformly rational), then using the weak factorisation theorem we could answer the question affirmatively. Unfortunately, so far this line of thought has not been completed.

The goal of this paper is to prove the following theorem, which brings us significantly closer to answering the original question of Gromov:
\begin{thm}\label{main_thm}
Let $X$ be a nonsingular variety, such that $X\times \C$ is rational. Then, $X\times \C^2$ is uniformly rational.
\end{thm}

As a consequence, we obtain a positive answer to the stably rational version of the problem:
\begin{cor}
Let $X$ be a nonsingular stably rational variety over $\C$. Then it is stably uniformly rational.
\begin{proof}
By assumption, there exists $n\geq 0$ such that $X\times \C^n$ is rational. Then according to Theorem \ref{main_thm}, $X\times \C^{n+2}$ is uniformly rational.
\end{proof}
\end{cor}

Some of the ingredients of our current approach are similar to what was applied in \cite{BaneckiUniformRetract} in the retract rational case; Lemma \ref{lem_generic} is a variant of \cite[Proposition 2.2]{BaneckiUniformRetract} and Lemma \ref{lem_denominator_red} is a variant of \cite[Lemma 2.4]{BaneckiUniformRetract}. Nonetheless, the current reasoning is more complicated than the one presented in \cite{BaneckiUniformRetract}. In particular, it depends also on a result of the paper \cite{MellaEquivBirEmbedIII2013} of Mella, which provides a condition under which there exists a Cremona automorphism of $\C^n$ transforming a given divisor into another. Here can be seen a connection between the problem of uniform rationality of rational varieties and certain properties of the Cremona group; for example in Section \ref{sec_Cremona} we propose a conjecture regarding the group, which would strengthen the conclusion of Theorem \ref{main_thm}.

\section{Blowups in special ideals}
In this section, we strengthen Observation \ref{obs_Bogomolov}, going slightly beyond the case of nonsingular centre. Before passing to the main result of this section, we introduce the following lemma, which is a variant of \cite[Proposition 2.2]{BaneckiUniformRetract}. It could be derived as a consequence of that result, nonetheless we believe that it is more informative to prove it directly:
\begin{lem}\label{lem_generic}
Let $(X,0)\subset(\C^n,0)$ be a nonsingular germ of a subvariety of $\C^n$ of dimension $m<n$. Let $I$ be a nonzero ideal of the ring $\OO_{X,0}$. After a generic linear change of coordinates, if $\pi$ denotes the projection to the first $m$ variables $\pi:(X,0)\rightarrow (\C^m,0)$, then the homomorphism induced by pre-composing with $\pi$ and the taking the quotient modulo $I$:
\begin{equation*}
    \pi^\ast:\OO_{\C^m,0}\rightarrow\OO_{X,0}/I
\end{equation*}
is surjective.
\begin{proof}
Since we are interested only in the germ of $X$ at $0$, we can assume that $X$ is an irreducible Zariski closed subset of $\C^n$. After the generic change of coordinates, we can assume that the derivative of $\pi$ at $0\in X$ is an isomorphism. In particular, $\pi$ induces a ring extension $\C[x_1,\dots,x_m]\subset \OO(X)$. From Noether normalisation lemma we can assume that the extension is finite.

Denote by $\mathfrak{n}$ the maximal ideal of functions vanishing at $0$ in $\OO(X)$ and by $\mathfrak{m}=\mathfrak{n}\cap \C[x_1,\dots,x_m]$ the maximal ideal of polynomials vanishing at $0$ in $\C[x_1,\dots,x_m]$. 

The local ring $\OO_{X,0}$ is the localisation of $\OO(X)$ in $\mathfrak{n}$, so we have the inclusion $\OO(X)\subset \OO_{X,0}$. Define $J:=I\cap \OO(X)$.

After counting dimensions (\cite[Lemma 2.3]{BaneckiUniformRetract}), we can assume that the only point of the zero set of $J$ lying in the preimage of $0$ through $\pi$ is $0$. This means that
\begin{equation*}
    \sqrt{\mathfrak{m}\OO(X)+J}=\mathfrak{n}.
\end{equation*}

As the derivative of $\pi$ at $0$ is an isomorphism, we have that the elements of $\mathfrak{m}$ generate the maximal ideal of the local ring $\OO_{X,0}$. This means that for every $f\in\mathfrak{n}$ there is $q\in\OO(X)\backslash\mathfrak{n}$ such that
\begin{equation*}
	qf\in \mathfrak{m}\OO(X)\subset \mathfrak{m}\OO(X) +J.
\end{equation*}
As the latter ideal is $\mathfrak{n}$-primary we get that 
\begin{equation}\label{ideal_equal_eq}
	\mathfrak{m}\OO(X) +J=\mathfrak{n}.
\end{equation}

Consider $\OO(X)$ as a $\C[x_1,\dots,x_m]$-module and define the quotient module $M:=\OO(X)/(I+\C[x_1,\dots,x_m])$. From \eqref{ideal_equal_eq} we deduce that 
\begin{equation*}
	M\subset M\mathfrak{m}.
\end{equation*}
Hence, it follows from Nakayama's lemma that $M_{\mathfrak{m}}=0$. This is equivalent to the conclusion of the lemma.
\end{proof}
\end{lem}

In the proof of the main result of this section we will need to work in specific blowup coordinates, and that is why we would like to recall precisely how changing the set of generators of a given ideal corresponds to a change of coordinates of the blown up variety. We are specifically concerned with the case of codimension two centres:
\begin{defiandobs}
Let $X$ be a variety and let $f,g\in \OO_{X,x_0}$ be two nonzero germs of regular functions at a point $x_0\in X$. 

Choose an open neighbourhood $U$ of $x_0$ on which both $f$ and $g$ are represented by global regular functions. We define the \emph{local blowup} $\Bl(f,g)$ of $X$ in $(f,g)$ as the Zariski closure in $U\times \mathbb{P}^1_\C$ of the set
\[
    \{(x,(f(x):g(x))) \in (U \setminus \Z(f,g)) \times \mathbb{P}^1_\C \}.
\]
The construction depends on the choice of $U$, but since we are only interested in properties of $\Bl(f,g)$ near the fibre over $x_0$ under the natural projection 
\(\sigma:\Bl(f,g) \to U\), the ambiguity will be irrelevant.

Let $(f',g')$ be another pair of functions generating the ideal $(f,g)$ in $\OO_{X,x_0}$. Then $\Bl(f,g)$ and $\Bl(f',g')$ are isomorphic. More specifically, shrinking $U$ if necessary, there exists a representation of $f'$ and $g'$ as a linear combination of $f$ and $g$ with coefficients in $\OO(U)$:
\begin{align*}
    f'=af+bg, \\
    g'=cf+dg.
\end{align*}
Then, the mapping 
\begin{align*}
    \varphi&:U\times \mathbb{P}_\C^1\rightarrow U\times \mathbb{P}_\C^1,\\
    \varphi(x,t)&:=\left(x,\begin{bmatrix}
        a(x) & b(x) \\
        c(x) & d(x)
    \end{bmatrix}(s:t)\right)
\end{align*}
restricts to an isomorphism between $\Bl(f,g)$ and $\Bl(f',g')$ (here the matrix is treated as an element of $\mathrm{PGL}(2,\C)$).
\end{defiandobs}

Using the introduced notation, we make the following simple observation:
\begin{obs}\label{obs_point}
Let $X$ be a variety, let $x_0$ be a point of $X$ and let $f,g\in \OO_{X,x_0}$. Assume that $\frac{g}{f}\not\in \OO_{X,x_0}$. Then, $(x_0,(0:1))\in \Bl(f,g)$.
\begin{proof}
If $g$ divides $f$ in $\OO_{X,x_0}$, then the mapping
\begin{equation*}
    \psi:U\rightarrow \Bl(f,g), \quad
    \psi(x):=\left(x,\left(\frac{f(x)}{g(x)}:1\right)\right)
\end{equation*}
defined on a neighbourhood $U$ of $x_0$ is an isomorphism onto $\sigma^{-1}(U)$, where $\sigma:\Bl(f,g)\rightarrow U$ is the projection. As $\frac{f}{g}$ is not invertible in $\OO_{X,x_0}$, we have that the value of the quotient at $x_0$ is equal to $0$, so $(x_0,(0:1))\in \Bl(f,g)$.

If $g$ does not divide $f$ in $\OO_{X,x_0}$, then the ideal $(f,g)$ is not principal in the local ring, so the rational morphism $\sigma^{-1}$ cannot be extended as a regular mapping near the point $x_0$. It follows from Zariski main theorem that $\sigma^{-1}(x_0)$ is of positive dimension, so it is equal to $\{x_0\}\times \mathbb{P}_\C^1$. In particular it contains the point $(x_0,(0:1))$.
\end{proof}
\end{obs}

We can now state and prove the main result of this section:
\begin{prop}\label{prop_blowup}
Let $f,g\in\OO_{\C^n,0}$ be two nonzero germs of regular functions at the origin in $\C^n$ for some $n$. Assume that $f$ vanishes at $0$, but its derivative does not, and that $f$ does not divide $g$ in the ring $\OO_{\C^n,0}$. Then, the germ
\begin{equation*}
    (\Bl(f,g),(0,(0:1)))
\end{equation*}
is isomorphic to $(\C^n,0)$.
\begin{proof}
First note that thanks to Observation \ref{obs_point}, the point $(0,(0:1))$ indeed belongs to $\Bl(f,g)$. 

If $g(0)\neq 0$ then the conclusion of the theorem is obvious, so we assume that $g(0)=0$. Let $(X,0)$ be the germ of the zero set of $f$ at $0$. After a generic change of coordinates, using Lemma \ref{lem_generic} we have that the projection $\pi:X\rightarrow \C^{n-1}$ satisfies
\begin{equation*}
    \pi^\ast:\OO_{\C^{n-1},0}\rightarrow \OO_{X,0}/(g^2)\text{ is a surjection}.
\end{equation*}
Moreover, since the coordinate change is generic we can assume that $\pdv{f}{x_n}\neq 0$. Let $h\in \OO_{\C^{n-1},0}$ be such that $\pi^\ast(h)\equiv x_n\mod g^2$. Then $h-x_n\in(f,g^2)$ in $\OO_{\C^n,0}$. This means that $h-x_n$ can be written as 
\begin{equation*}
    h-x_n=af+bg^2
\end{equation*}
for some $a,b\in \OO_{\C^n,0}$. We must have that $a(0)\neq 0$, for otherwise $0=\pdv{h-x_n}{x_n}=-1$. It follows that $f\in (h-x_n,g)$, so $(h-x_n,g)=(f,g)$ in $\OO_{\C^n,0}$. The transition matrix between the two sets of generators is given by
\begin{equation*}
    M=\begin{bmatrix}
        a & bg \\
        0 & 1
    \end{bmatrix}.
\end{equation*}
After evaluating the functions at the point $0$ the matrix satisfies $M(0:1)=(0:1)$, so the germs $(\Bl(f,g),(0,(0:1)))$ and $(\Bl(h-x_n,g),(0,(0:1)))$ are isomorphic. Thus, after replacing $f$ by $h-x_n$, without loss of generality we may assume that $f$ is of the form $h-x_n$ for some $h\in\OO_{\C^{n-1},0}$. 

The coordinate change
\begin{equation*}
    \varphi:(\C^n,0)\rightarrow (\C^n,0),\quad
    \varphi(x):=(x_1,\dots,x_{n-1},h-x_n)
\end{equation*}
is an involution, hence a biregular automorphism of the germ. When applied to the coordinates on $\C^n$, $\varphi$ transforms $f$ into $x_n$, so now without loss of generality we can take $f=x_n$.

We can write down $g=g_1+x_ng_2$ for some $g_1\in \OO_{\C^{n-1},0}, g_2\in\OO_{\C^n,0}$. It follows easily that $(x_n,g_1)=(x_n,g)$ in $\OO_{\C^n,0}$ and the transition matrix between the two sets of generators is given by
\begin{equation*}
    N=\begin{bmatrix}
        1 &0 \\
        g_2 & 1
    \end{bmatrix}.
\end{equation*}
Again, we have that $N(0:1)=(0:1)$, so substituting $g_1$ for $g$ we assume that $g\in \OO_{\C^{n-1},0}$.

Finally, we can explicitly define an isomorphism $\psi:(\Bl(x_n,g),(0:1))\rightarrow (\C^n,0)$ by
\begin{equation*}
    \psi(x,(t:s)):=\left(x_1,\dots,x_{n-1},\frac{s}{t}\right),
\end{equation*}
The inverse of the morphism is then given by
\begin{align*}
    \psi^{-1}&:(\C^n,0)\rightarrow (\Bl(x_n,g),(0:1)), \\
    \psi^{-1}(y_1,\dots,y_n)&:=(y_1,\dots,y_{n-1},y_n g(y_1,\dots,y_{n-1}),(y_n:1)).
\end{align*}
\end{proof}
\end{prop}

\section{Reduction to a Cremona equivalence problem}
In this section we will give a certain condition regarding existence of a particular element of the Cremona group of birational transformations of the affine space $\C^n$, which would imply Theorem \ref{main_thm}. Before that, we need a few preliminary results. The following lemma appeared earlier as \cite[Lemma 3.2]{BaneckiRelativeSW}, but since it was stated and proven only over the real numbers, we repeat the proof here.
\begin{lem}\label{lem_denominator_red}
Let $(X,x_0)$ be a germ of a variety. Let $f:X\times \C\dashrightarrow \C$ be a rational function, which is defined at the generic point of $X\times \{0\}$. Assume that the restriction $f\vert_{X\times\{0\}}$ extends to a regular germ at $x_0$. Let $Q$ be a denominator of $f$, i.e. a nonzero regular germ $Q\in\OO_{X\times \C,(x_0,0)}$, such that the product $Qf$ is regular at $(x_0,0)$. Then, the rational function
\begin{equation*}
    g:X\times \C\dashrightarrow \C,\quad
    g(x,t):=f(x,Q(x,0)t)
\end{equation*}
is regular at $(x_0,0)$.
\begin{proof}
Denote by $F\in \OO_{X,x_0}$ the restriction $f\vert_{X\times \{0\}}$ as a function of $x$, which by assumption is regular at $x_0$. After substituting $f'(x,t):=f(x,t)-F(x)$ without loss of generality we can assume that $f$ vanishes on $X\times \{0\}$. This means that the regular function $Qf\in\OO_{X\times \C,(x_0,0)}$ can be written as $Qf=tR$ for some $R\in \OO_{X\times \C,(x_0,0)}$, where $t$ is the coordinate on the factor $\C$. Similarly, $Q$ can be written as $Q(x,t)=Q(x,0)+tS(x,y)$ for some $S\in\OO_{X\times \C,(x_0,0)}$. Substituting everything into the equation defining $g$ we find that
\begin{multline*}
    f(x,Q(x,0)t)= \frac{Q(x,Q(x,0)t)f(x,Q(x,0)t)}{Q(x,Q(x,0)t)}=\\
    =\frac{Q(x,0)tR(x,Q(x,0)t)}{Q(x,0)+ Q(x,0)tS(x,Q(x,0)t)} 
    = \frac{tR(x,Q(x,0)t)}{1+ tS(x,Q(x,0)t)}.
\end{multline*}
The above denominator does not vanish at $(x_0,0)$, hence the fraction gives a regular representation of $g$ on a Zariski neighbourhood of that point.
\end{proof}
\end{lem}

We will also make use of the following two simple facts:
\begin{obs}\label{obs_etale}
Let $(X,x_0)$ and $(Y,y_0)$ be two nonsingular germs of varieties of the same dimension. Let $\eta:(X,x_0)\rightarrow (Y,y_0)$ be a germ of a regular mapping, which is a birational isomorphism. Assume that the derivative of $\eta$ at $x_0$ is an isomorphism. Then $\eta$ is an isomorphism between the germs $(X,x_0)$ and $(Y,y_0)$.
\begin{proof}
Let $\mu:Y\dashrightarrow X$ be the rational inverse of $\eta$. The composition $\mu\circ\eta$ is equal to the identity, so in particular it extends to a regular mapping at $x_0$. It follows from \cite[Observation 2.1]{BaneckiUniformRetract} that $\mu$ extends to a regular mapping at $\eta(x_0)=y_0$, which gives the local regular inverse of $\eta$. 
\end{proof}
\end{obs}

\begin{obs}\label{obs_iso_at_generic}
Let $\psi:X\dashrightarrow Y$ be a birational isomorphism between two varieties $X$ and $Y$. Let $Z$ and $W$ be subvarieties of codimension one of $X$ and $Y$ respectively, such that $\psi$ is defined at the generic point of $Z$ and when restricted to $Z$ is a birational isomorphism onto $W$. Then $\psi$ is an isomorphism at the generic point of $Z$.
\begin{proof}
Without loss of generality we can take projective closures of $X$ and $Y$ to assume that they are projective. Now, as $X$ is projective the indeterminancy locus of $\psi^{-1}:Y\dashrightarrow X$ is of codimension at least two. In particular $\psi^{-1}$ is defined at the generic point of $W$, so it gives a local inverse of $\psi$ at the generic point of $Z$.
\end{proof}
\end{obs}

Our goal now is to prove the following theorem, which gives a sufficient condition for a variety to be uniformly rational. Later, we will verify that the condition is satisfied in the context of Theorem \ref{main_thm}.

In the theorem and its proof, we will work locally with germs of varieties and regular mappings. In the diagrams, dashed arrows signify rational mappings, while solid ones denote germs of regular mappings. It is always clear in a neighbourhood of which point the germs are defined, so we omit this information in the pictures.

\begin{thm}\label{thm_reduction}
Let $(X,x_0)$ be a nonsingular germ of a variety of dimension $n$. Assume that there exists an embedding $\varphi:(X,x_0)\hookrightarrow (\C^{n+1},0)$, and a birational isomorphism $\psi:X\times \C\dashrightarrow \C^{n+1}$, which is defined at the generic point of $X\times \{0\}$, such that the diagram
\begin{center}
\begin{tikzcd}
X \arrow[r, "\varphi", hook] \arrow[d, hook] & \C^{n+1} \\
X\times \C \arrow[ru, "\psi", dashed]  &         
\end{tikzcd}
\end{center}
is commutative. Here the arrow $X\hookrightarrow X\times \C$ denotes the canonical embedding identifying $X$ with $X\times \{0\}$.

Then, the germ $(X\times \C,(x_0,0))$ is isomorphic to $(\C^{n+1},0)$.
\begin{proof}
\begin{claim}\label{claim1}
Without loss of generality we can assume that $\psi$ is defined at $(x_0,0)$.
\begin{proof}
Let $Q\in \OO_{X\times \C,(x_0,0)}$ be a common denominator of all the coordinates of $\psi$. As $\psi$ is defined at the generic point of $X\times \{0\}$, we can assume that $Q$ does not vanish identically on $X\times \{0\}$. According to Lemma \ref{lem_denominator_red}, if we define
\begin{gather*}
    \xi:(X\times \C,(x_0,0)) \rightarrow (X\times \C,(x_0,0)), \\
    \xi(x,t):=(x,Q(x,0)t),
\end{gather*}
then we have that $\psi\circ \xi$ extends to a regular germ at $(x_0,0)$. Note that $\xi$ is a birational isomorphism with the inverse given by $(x,t)\mapsto (x,\frac{t}{Q(x,0)})$. Note also that $\xi$ satisfies $\xi(x,0)=(x,0)$ for $x\in X$, so the following diagram is commutative:
\begin{center}
\begin{tikzcd}
X \arrow[r, "\varphi"] \arrow[d, hook] & \C^{n+1} \\
X\times \C \arrow[d, "\xi"]  \arrow[ru]          &          \\
X\times \C \arrow[ruu, "\psi"', dotted] &         
\end{tikzcd}
\end{center}

Now we can just take the composition $\psi\circ \xi$ instead of $\psi$. This finishes the proof of the claim.
\end{proof}
\end{claim}

Let now $f\in\OO_{\C^n,0}$ be a local generator of the ideal of germs vanishing on the image of $(X,x_0)$ through $\varphi$. The germ is nonsingular, so the derivative of $f$ does not vanish at the point $0$. The pullback $f\circ \psi$ vanishes on $X$, so it can be written as:
\begin{equation*}
    f\circ \psi=tF,
\end{equation*}
where $F\in \OO_{X\times \C,(x_0,0)}$ and $t$ is the coordinate on the factor $\C$. According to Observation \ref{obs_iso_at_generic}, $\psi$ is an isomorphism at the generic point of $X\times \{0\}$, so we have that the derivative of $f\circ \psi$ does not vanish identically on $X\times \{0\}$, hence neither does $F$.

\begin{claim}\label{claim2}
Without loss of generality we can assume that there exists $H\in\OO_{X,x_0}$ such that
\begin{enumerate}
    \item $F$ can be written as $F(x,t)=H(x)u(x,t)$, where $u\in\OO_{X\times \C,(x_0,0)}^\ast$ is a unit,
    \item $\psi(x,t)=\varphi(x)$ for $x\in \Z(H),t\in \C$.
\end{enumerate}
\begin{proof}
To begin with, note that $F$ can be written as 
\begin{equation}\label{thm_reduction_eq0}
    F(x,t)=F(x,0)+tG(x,t)
\end{equation}
for some $G\in\OO_{X\times \C,(x_0,0)}$.

Similarly to what was done in the proof of Claim \ref{claim1}, define 
\begin{gather*}
    \xi:(X\times \C,(x_0,0))\rightarrow (X\times \C,(x_0,0)), \\
    \xi(x,t):=(x,F(x,0)t).
\end{gather*}
Then, using \eqref{thm_reduction_eq0}, we can factor the composition $f\circ\psi\circ \xi$ into the following product:
\begin{multline*}
    f\circ\psi\circ \xi(x,t)=tF(x,0)F(x,tF(x,0))=\\
    =tF(x,0)(F(x,0)+tF(x,0)G(x,tF(x,0)))=tF^2(x,0)(1+tG(x,tF(x,0))).
\end{multline*}

Now, as it was done in the proof of Claim \ref{claim1}, we will substitute the composition $\psi\circ \xi$ in place of $\psi$. One can then take 
\begin{align*}
    H(x)&:=F^2(x,0),\\
    u(x,t)&:=(1+tG(x,F(x,0)t)).
\end{align*}
We do have
\begin{equation*}
    f\circ \psi\circ \xi(x,t)=tH(x)u(x,t).
\end{equation*}
Furthermore, note that for $x\in \Z(H)$ we have 
\begin{equation*}
    \psi\circ\xi(x,t)=\psi(x,F(x,0)t)=\psi(x,0)=\varphi(x).
\end{equation*}
This finishes the proof of Claim \ref{claim2}.
\end{proof}
\end{claim}

\begin{claim}\label{claim3}
There exists $\overline{H}\in\OO_{\C^{n+1},0}$, such that $\overline{H}\circ \psi \sim H$ in $\OO_{X\times \C,(x_0,0)}$ (here tilde signifies that the two elements differ by multiplication by a unit).
\begin{proof}
Write down $H$ as a product of primes in $\OO_{X,x_0}$:
\begin{equation*}
    H=H_1\dots H_k.
\end{equation*}

As $\varphi$ is an isomorphism onto its image, there do exist $\overline{H}_1,\dots,\overline{H}_k\in\OO_{\C^{n+1},0}$ which satisfy
\begin{equation*}
    \overline{H}_i\circ \varphi =H_i\text{ for } 1\leq i \leq k,
\end{equation*}
as elements of $\OO_{X,x_0}$. Note that, by point $2$ of Claim \ref{claim2}, for $x\in \Z(H_i)$ we have that
\begin{equation*}
    \overline{H}_i\circ \psi(x,t)=\overline{H}_i\circ\varphi(x)=H_i(x)=0.
\end{equation*}

As $H_i$ is a prime element of $\OO_{X,x_0}$, it is also prime when treated as an element of $\OO_{X\times \C,(x_0,0)}$. It follows that it divides $\overline{H}_i\circ \psi$.

On the other hand, the quotient $\frac{\overline{H}_i\circ \psi }{H_i}$ is constantly equal to one, when restricted to $X\times \{0\}$. It hence follows that we must have $H_i\sim \overline{H}_i\circ \psi$ in $\OO_{X\times \C,(x_0,0)}$. Now, we can take $\overline{H}:=\overline{H}_1\dots \overline{H}_k$. This finishes the proof of Claim \ref{claim3}.
\end{proof}
\end{claim}

We can now finish the proof of Theorem \ref{thm_reduction}. Consider the ideal $(f,\overline{H})$ in $\OO_{\C^{n+1},0}$. Note that $\overline{H}$ does not vanish identically on the zero set of $f$, as this would imply that $H$ is constantly equal to zero. Using the notation of Proposition \ref{prop_blowup}, consider the germ of the blowup of $\C^{n+1}$ in the ideal $(f,\overline{H})$ together with the projection $\sigma:(\Bl(f,\overline{H}),(0,(0:1)))\rightarrow (\C^{n+1},0)$. Note that $\overline{H}$ was chosen in such a way that
\begin{equation}\label{thm_reduction_eq1}
    \frac{f\circ \psi}{\overline{H}\circ \psi}=tu',
\end{equation}
where $u'\in \OO_{X\times \C,(x_0,0)}^\ast$ is a unit. This means that the rational mapping $\eta:X\times \C\dashrightarrow \Bl(f,\overline{H})$ making the diagram below commutative is defined at $(x_0,0)$, and maps the point to $(0,(0:1))$:
\begin{center}
\begin{tikzcd}
                                                & {\Bl(f,\overline{H})} \arrow[d, "\sigma"] \\
X\times \C \arrow[r, "\psi"] \arrow[ru, "\eta"] & \C^{n+1}                              
\end{tikzcd}
\end{center}

We claim that $\eta$ is an isomorphism at $(x_0,0)$, in the view of Proposition \ref{prop_blowup} this will show that 
\begin{equation*}
    (X\times \C,(x_0,0))\cong (\Bl(f,\overline{H}),(0,(0:1)))\cong (\C^{n+1},0). 
\end{equation*}

Since $\eta$ is birational, according to Observation \ref{obs_etale}, it suffices to check that its derivative at the point $(x_0,0)$ is an isomorphism. Let $v\in T_{(x_0,0)}(X\times \C)$ belong to the kernel of the derivative. Equation \eqref{thm_reduction_eq1} shows that the element $tu'\in \OO_{X\times \C,(x_0,0)}$ lies in the image of the induced homomorphism
\begin{align*}
    \eta^\ast &:\OO_{\Bl(f,\overline{H}),(0,(0:1))}\rightarrow \OO_{X\times \C,(x_0,0)},\\
    \eta^\ast(g)&:=g\circ \eta.
\end{align*}
In particular, $v$ must lie in the kernel of the derivative of $tu'$ at $(x_0,0)$, which is precisely the tangent space to $X\times \{0\}$.

However, by assumption the restriction $\sigma\circ\eta\vert_{X\times\{0\}}=\psi\vert_{X\times \{0\}}=\varphi$ is an immersion, so $v=0$. This finishes the proof.
\end{proof}
\end{thm}

We will now provide a condition, under which the assumptions of Theorem \ref{thm_reduction} are automatically satisfied. To do that we recall the following notion:
\begin{defi}
Let $X$ and $Y$ be two subvarieties of codimension one embedded in $\C^n$. We say that they are \emph{Cremona equivalent}, if there exists a birational automorphism $\alpha:\C^n\rightarrow \C^n$, which is defined at the generic point of $X$, and when restricted to $X$ is a birational isomorphism onto $Y$.
\end{defi}

\begin{prop}\label{prop_cremona=>ur}
Let $(X,x_0)$ be a nonsingular germ of a rational variety of dimension $n$. Assume that there exists an embedding of the germ $\varphi:(X,x_0)\hookrightarrow \C^{n+1}$, such that the image $\varphi(X)$ is Cremona equivalent to a hyperplane. Then, the germ $(X\times \C,(x_0,0))$ is isomorphic to $(\C^{n+1},0)$. 
\begin{proof}
Let $\alpha$ be the Cremona automorphism transforming $\varphi(X)$ into the hyperplane, and denote by $\psi_0$ the composition $\alpha\circ \varphi$. Then, the assumptions are presented in the following diagram:
\begin{center}
\begin{tikzcd}
                                                     & \C^{n+1} \arrow[d, "\alpha", dashed] \\
                                                     & \C^{n+1}                             \\
X \arrow[ruu, "\varphi", hook] \arrow[r, "\psi_0", dashed] & \C^n \arrow[u, hook]                
\end{tikzcd}
\end{center}

If we now define 
\begin{equation*}
    \psi_1:X\times \C\dashrightarrow \C^{n+1}, \quad
    \psi_1(x,t):=(\psi_0(x),t)
\end{equation*}
we have that $\psi_1$ is a birational isomorphism and the following diagram is commutative:
\begin{center}
\begin{tikzcd}
X\times \C \arrow[r, "\psi_1", dashed]                & \C^{n+1}             \\
X \arrow[r, "\psi_0", dashed] \arrow[u, hook] & \C^n \arrow[u, hook]
\end{tikzcd}
\end{center}

Together, from the two diagrams we have commutativity of the following:
\begin{center}
\begin{tikzcd}
X\times \C \arrow[r, "\psi_1", dashed] & \C^{n+1}  \\
X \arrow[r, "\varphi", hook] \arrow[u, hook] & \C^{n+1} \arrow[u, "\alpha", dashed]                         
\end{tikzcd}
\end{center}

Hence, if we define $\psi:=\alpha^{-1}\circ \psi_1$ we have that the assumptions of Theorem \ref{thm_reduction} are satisfied.
\end{proof}
\end{prop}

\section{On Cremona equivalences}\label{sec_Cremona}
We will now deduce Theorem \ref{main_thm} from Proposition \ref{prop_cremona=>ur}. Before we do that, let us state the following question, which is natural to consider here:
\begin{que}\label{que1}
Let $X$ be a rational variety of dimension $n$ and let $x_0\in X$ be its nonsingular point. Does there always exist an embedding $\varphi:(X,x_0)\hookrightarrow \C^{n+1}$ of the germ, such that the image $\varphi(X)$ is Cremona equivalent to a hyperplane?
\end{que}

We do not know the answer to the question, but we conjecture that it is positive. In view of Proposition \ref{prop_cremona=>ur}, it would follow that for such a germ, the product $(X\times \C,(x_0,0))$ is isomorphic to $(\C^{n+1},0)$, which would strengthen Theorem \ref{main_thm}. On the other hand, note that if all nonsingular rational varieties are uniformly rational, then the answer to this question is positive; in such case $(X,x_0)$ can be embedded as a hypersurface in $\C^{n+1}$ and the Cremona transformation can just be taken to be the identity.

Question \ref{que1} is a special case of the following more general one:
\begin{que}\label{que2}
Let $X$ be a variety of dimension $n$ and let $x_0\in X$ be its nonsingular point. Let $\mu:X\dashrightarrow \C^{n+1}$ be a fixed birational embedding of $X$ into the affine space. Does there exist an embedding $\varphi:(X,x_0)\hookrightarrow \C^{n+1}$, and a Cremona automorphism $\alpha:\C^{n+1}\dashrightarrow \C^{n+1}$, such that the composition $\alpha\circ \varphi:X\dashrightarrow \C^{n+1}$ is well defined as a rational mapping, and equal to $\mu$? 
\end{que}
While conjecturing a positive answer to Question \ref{que2} may be overly optimistic, either a positive or negative resolution of Questions \ref{que1} and \ref{que2} would improve our insight into Question \ref{main_que}.

Fortunately, using results of \cite{MellaEquivBirEmbedIII2013}, we can provide an affirmative answer to the Question \ref{que1} in a particular case, which is sufficient to tackle Theorem \ref{main_thm}. Let us restate the main result of that paper:
\begin{thm}\label{thm_cones}
Let $X\subset \C^{n+1}$ and $Y\subset\C^{n+1}$ be two subvarieties of $\C^{n+1}$ of dimension $n$, which are birational. Then the varieties $X\times \C\subset \C^{n+2}$ and $Y\times \C\subset \C^{n+2}$ are Cremona equivalent.
\begin{proof}
Using the terminology of \cite{MellaEquivBirEmbedIII2013}, the varieties $X\times \C\subset \C^{n+2}$ and $Y\times \C\subset \C^{n+2}$ are cones over $X$ and $Y$, respectively, with vertices lying at infinity in the corresponding projective closures. As $X$ and $Y$ are birational, according to \cite[Theorem 2]{MellaEquivBirEmbedIII2013} their cones are Cremona equivalent.
\end{proof}
\end{thm}
\begin{thm}
Let $(X,x_0)$ be a nonsingular germ of a variety of dimension $n$, such that $X\times \C$ is rational. Then, there exists an embedding 
\begin{equation*}
    \varphi:(X\times \C,(x_0,0))\hookrightarrow \C^{n+2},
\end{equation*}
such that $\varphi(X\times \C)$ is Cremona equivalent to a hyperplane.
\begin{proof}
According to \cite[Proposition 2.1]{BodnarHauserSchichoVillamayor2008}, every nonsingular variety can be locally embedded as a hypersurface, so there exists an embedding $\varphi:(X,x_0)\hookrightarrow \C^{n+1}$. Consider the natural extension:
\begin{equation*}
    \Phi:(X\times \C,(x_0,0)) \hookrightarrow \C^{n+2}, \quad
    \Phi(x,t):=(\varphi(x),t).
\end{equation*}
We claim that the image of $X\times \C$ through $\Phi$ is Cremona equivalent to a hyperplane.

To show that, fix a birational isomorphism $\psi:X\times \C \dashrightarrow \C^{n+1}$. After a translation in the second variable, we can assume that $\psi$ is defined at the generic point of $X\times \{0\}$ and that the restriction $\psi\vert_{X\times \{0\}}$ is a birational isomorphism onto its image. Similarly as before, consider the birational extension of $\psi$ defined by
\begin{equation*}
    \Psi:X\times \C^2 \dashrightarrow \C^{n+2},\quad
    \Psi(x,t):=(\psi(x),t).
\end{equation*}

Now, thanks to Theorem \ref{thm_cones}, the varieties $\Phi(X\times \C)$ and $\Psi(X\times\{0\}\times \C)$ are Cremona equivalent, so now it is enough to show that $\Psi(X\times \{0\}\times\C)$ is Cremona equivalent to a hyperplane.

Let $\tau:X\times \C^2\rightarrow X\times \C^2$ be the automorphism permuting the last two coordinates of the product. We define a Cremona automorphism $\alpha:\C^{n+2}\dashrightarrow \C^{n+2}$ by the equation
\begin{equation*}
    \alpha:=\Psi\circ \tau \circ \Psi^{-1}.
\end{equation*}
It is then clear that $\alpha$ is a Cremona equivalence between the varieties $\Psi(X\times \{0\}\times \C)$ and $\Psi(X\times \C \times\{0\})$. The latter one by definition is a hyperplane in $\C^{n+2}$. This finishes the proof.
\end{proof}
\end{thm}

Together with Proposition \ref{prop_cremona=>ur}, this proves Theorem \ref{main_thm}.

\printbibliography
\end{document}